\documentclass[preprint,12pt]{elsarticle}

\usepackage{amssymb}
\usepackage{amscd,amsmath,amsfonts,graphicx}
\usepackage{amsthm}

\newcommand{\R}{\mathbb{R}}

\newcommand{\E}{{\bf E}}
\newcommand{\I}{\mathbb{I}}

\begin{document}
\baselineskip=22pt \centerline{\large \bf Necessary and Sufficient
Conditions for Unique Solution} \centerline{\large \bf to Functional
Equations of Poincar\'e Type}

\vspace{1cm} \centerline{Chin-Yuan Hu $^a$ and Gwo Dong Lin $^b$}
\centerline{$^a$ National Changhua University of Education, Taiwan,
and} \centerline{$^b$ Hwa-Kang Xing-Ye Foundation\  and  Academia
Sinica, Taiwan} \vspace{1cm} \noindent {\bf Abstract.}
Distributional equation is an important tool in the characterization
theory
 because many characteristic properties of
distributions can be transferred to such equations.  Using a novel
and natural approach, we retreat a remarkable distributional
equation whose corresponding functional equation in terms of
Laplace--Stieltjes transform is of the Poincar\'e type. The
necessary and sufficient conditions for the equation to have a {\it
unique} distributional solution with finite {\it variance} are
provided. This complements the previous results which involve at
most the {\it mean} of the distributional solution. Besides, more
general distributional (or functional) equations are investigated as
well.
\vspace{0.1cm}\\
\hrule
{\bf Short title: Poincar\'e Type Functional Equations}\\
{\bf Postal addresses:} Chin-Yuan Hu, Department of Business
Education, National Changhua University of Education, Changhua
50058, Taiwan. (E-mail: buhuua@gmail.com)\\
Gwo Dong Lin, (1) Social and Data Science Research Center, Hwa-Kang
Xing-Ye Foundation, Taipei 10659, Taiwan, and (2) Institute of
Statistical Science, Academia Sinica, Taipei 11529, Taiwan. (E-mail:
 gdlin@stat.sinica.edu.tw)
\newpage

\noindent{\bf 1. Introduction}
\newcommand{\bin}[2]{
   \left(
     \begin{array}{@{}c@{}}
         #1  \\  #2
     \end{array}
   \right)          }

 One useful method to characterize probability distributions is through suitable distributional equations
 (see, e.g., Ramachandran and Lau 1991, Rao et al.\,1994,  Lin 1994, Hu and Lin 2001, 2018, and the references therein).
 In this paper, we will retreat a remarkable distributional equation described below.

Let $X$ and $T$ be two nonnegative random variables having
distributions $F$ and $F_T$, respectively, denoted $X\sim F_X=F$,
$T\sim F_T.$ Let $\{X_i\}_{i=1}^{\infty}$ be a sequence of
independent and identically distributed (i.i.d.) random variables
having distribution $F$ on $\R_+\equiv[0,\infty),$ and let
$\{T_i\}_{i=1}^{\infty}$ be another sequence of i.i.d. random
variables having distribution $F_T.$ Moreover, let $N$ be a random
variable taking values in $\mathbb{N}_0\equiv \{0,1,2,\dots\},$ and
assume that all the random variables $X, X_i, T, T_i, N$ are
independent. For given $T$ and $N,$ we will investigate the
distributional equation
\begin{eqnarray}
X\stackrel{\rm d}{=}\sum_{i=1}^NT_iX_i,
\end{eqnarray}
where `\,$\stackrel{\rm d}{=}$\,' means equality in distribution and
the summation is zero if $N$ takes value 0. For applications of Eq
(1) in various fields, see, e.g., the survey paper by Liu (1997).

Let $P_N$ denote the probability generating function (pgf) of $N$
and let $\hat{F}$ be the Laplace--Stieltjes transform of $X\sim F$;
namely, $P_N(t)=\E[t^N]=\sum_{k=0}^{\infty}\Pr(N=k)t^k,\ t\in[0,1],$
where $0^0\equiv 1,$ and $\hat{F}(s)=\E[\exp(-sX)],$ $s\ge 0.$ Then
the distributional equation (1) can be further transferred to the
following functional equation in terms of $\hat{F}, F_T$ and $P_N$:
\begin{eqnarray}
\hat{F}(s)=P_N\bigg(\int_0^{\infty}\hat{F}(ts)dF_T(t)\bigg)=P_N\big(\E[\exp(-sTX)]\big),\ \ s\ge 0.
\end{eqnarray}
When $F_T$ is a degenerate distribution at $p\in(0,1),$ namely,
$\Pr(T=p)=1,$  Eq (1) reduces to\  $X\stackrel{\rm
d}{=}\sum_{i=1}^NpX_i,$ and Eq (2) is exactly the Poincar\'e
functional equation
\begin{eqnarray}\hat{F}(s)=P_N\big(\hat{F}(ps)\big),\ \ s\ge 0,\end{eqnarray}
which arises in the Galton--Watson processes (Poincar\'e 1886,
1890). So we call the general Eq (2) a functional equation of the
Poincar\'e type.

It is seen that once Eq (1) or Eq (2) has a solution $X\sim F,$ each
constant multiplication of $X$ also plays a solution to Eq (1).
However, the solution might be unique, provided we fix the mean of
the distributions. In this sense, Eq (1) or Eq (2) becomes a
characteristic property of the distributional solution.  A typical
example is the classical characterization of the exponential
distribution through Eq (3), where we can take $N$ obeying the
geometric distribution: $\Pr(N=n)=p(1-p)^{n-1},\ n\ge 1;$ see, e.g.,
Azlarov and Volodin (1986), p.\,79.

The properties of the solutions $X\sim F$ heavily depend on those of
the given $T$ and $N.$ Some results about Eqs (1) and (2) are
available in the literature. For example, denote the counting number
$\tilde{N}=\sum_{i=1}^N\I_{\{T_i>0\}},$ where $\I_A$ is the
indicator function of the set $A.$ Then for given $N$ and $T\sim
F_T$ with the conditions $\Pr(\tilde{N}=0\ \hbox{or}\ 1)<1$ and
$\Pr(T=0\ \hbox{or}\ 1)<1$ (which are used to exclude some trivial
cases),
the following results hold (see Liu 2002, Theorem 1.1,  and the references therein):\\
 (i) Eq (1) (or Eq (2)) has a
solution $0\le X\sim F$ iff the random variables $N$ and $T$
together satisfy the conditions
\begin{eqnarray}
\Pr(T>0)\E[N]>1,\ \ \E[N]\E[T^{\alpha}]=1\ \ \hbox{and}\ \ \
\E[T^{\alpha}\log T]\le 0\ \ \ \hbox{for\ some}\
\alpha\in(0,1];
\end{eqnarray}
 (ii) Eq (1) (or Eq (2)) has a solution $0\le X\sim F$  with finite
{\it mean} iff the random variables $N$ and $T$ together satisfy the
conditions
\begin{eqnarray}
\Pr(T>0)\E[N]>1,\ \ \E[N]\E[T]=1,\ \ \E[N\log^+N]<\infty\ \ \ \hbox{and}\ \ \
\E[T\log T]< 0,
\end{eqnarray}
 where $\log^+ x=\log x$ if $x\ge 1$ and
$\log^+x=0$ otherwise, and $0\log 0\equiv 0.$

One of the main purposes in this paper is to find the necessary and
sufficient conditions for which Eq (2) (or Eq (1)) has a {\it
unique} solution $F$ (on $\mathbb{R}_+$) with a fixed mean and
finite {\it variance}, and hence it can be used to characterize
distributions. This complements the above results (i) and (ii) which
involve at most the {\it mean} of the distributional solution. Our
approach is different from the previous ones and is somehow more
natural. Moreover, some general cases are also investigated. The
main results are stated in the next section, while their proofs  are
given in Section 4.  The needed lemmas are provided in Section 3.
Finally, we have some discussions in Section 5.
\bigskip\\
 \noindent{\bf 2. Main results}

 We start with the simplest case Eq (1) (or Eq (2)). More complicated cases will follow.

\noindent{\bf Theorem 1.} Let $0\le X\sim F$ with Laplace--Stieltjes
transform $\hat{F}$ and $\mu=\E[X]\in(0,\infty).$ Let $T\ge 0$ and
$N\ge 0$ be two given random variables, where $N$ takes values in
$\mathbb{N}_0$ and has pgf $P_N.$  Then for given $\mu,$ the random
variables $N$ and $T$ together satisfy the conditions
\begin{eqnarray} \E[N]\E[T]=1,\ \ 0<\E[T^2]<\E[T]<1,\ \ \hbox{and}\ \
\E[N^2]<\infty \end{eqnarray} iff the functional equation (2) has
exactly one solution $F$ with mean $\mu$ and a finite variance.
Moreover, the variance is of the form
\begin{eqnarray}\hbox{Var}(X)=\frac{(\E[T])^2
\hbox{Var}(N)+\E[N] \hbox{Var}(T)}{1-\E[N]\E[T^2]} \cdot \mu^2
\end{eqnarray}
with $\E[N]=1/{\E[T]}.$
\medskip\\
\indent Unlike the previous results (i) and (ii) (which assume some
initial conditions on $\tilde{N}$ and $T$ to exclude the trivial
cases), we don't assume explicitly any initial condition in Theorem
1. But each of  the sufficiency and necessity parts does imply
implicitly the following: $\Pr(N=1)<1$ or $\Pr(T=1)<1.$ To see this,
if on the contrary $\Pr(N=1)=\Pr(T=1)=1,$ then the second condition
in (6) fails to hold. Moreover, in this case,
 Eqs (1) and (2) reduce to the identities
$X\stackrel{\rm d}{=}X_1$ and $\hat{F}(s)=\hat{F}(s),\ s\ge 0,$
respectively, so the solution to Eq (2) is not unique, a
contradiction to the assumption in the sufficiency part.

 The following
result is about a characterization of degenerate distributions.

 \noindent{\bf Corollary 1.} Under the setting of Theorem 1, the functional equation (2) has
exactly one solution $F$ degenerate at mean $\mu=\E[X]$ iff the
random variables $N$ and $T$ are degenerate at $\E[N]$ and $\E[T],$
respectively, and $\E[T]=(\E[N])^{-1}\in(0,1);$ precisely,
$\Pr(N=n_0)=1$ for some integer $n_0\ge 2$ and $\Pr(T=1/n_0)=1.$

When $\Pr(T=p)=1$ for some $p\in(0,1)$ in Theorem 1, we are able to
rewrite Hu and Cheng's (2012) Theorem 1 with $\alpha=1$ and
Corollary 1 as follows.

 \noindent{\bf Corollary 2.} Let $p\in(0,1)$ and $\mu\in(0,\infty)$ be two constants.
 Let $N\ge 0$ be a random variable taking values in $\mathbb{N}_0$ and
 let $0\le X\sim F$ with mean $\mu$ and Laplace--Stieltjes transform $\hat{F}.$
 Then for given $\mu,$ the random variable $N$ satisfies the conditions
 \[\E[N]=1/p\ \ \hbox{and}\ \ \E[N^2]<\infty\]
 iff the Poincar\'e functional equation (3)  has exactly one
solution $F$ with mean $\mu$ and a finite variance. Moreover, the
variance is  equal to
\[ \hbox{Var}(X)=\frac{p^2 \hbox{Var}(N)}{1-p} \cdot \mu^2.\]

It is seen that the set of  conditions (5) is stronger than (4),
while (6) is stronger than both (4) and (5), as seen below. This in
turn implies that the solution $F$ in Theorem 1 belongs to the
classes
 of the previous solutions to Eq (2) under conditions (4) or (5).\\
\noindent{\bf Proposition 1.} Suppose that $N$ and $T$ are two
nonnegative random variables satisfying $\E[N]\E[T]=1$ and
$0<\E[T^2]<\E[T]<1$. Then $\Pr(T>0)\E[N]>1$ {and} $\E[T\log T]<0.$
\medskip\\
\indent If, in addition to (1),  assume that $N\ge m,$ where $m\ge
1$ is an integer, then we can split the RHS of (1) into two parts:
\begin{eqnarray}
X\stackrel{\rm d}{=}\sum_{i=1}^mT_iX_i+\sum_{j=1}^{N^*}T_{j+m}X_{j+m},
\end{eqnarray}
where $N^*=N-m\ge 0.$ This is equivalent to study the functional
equation
\begin{eqnarray}
\hat{F}(s)&=&\bigg(\int_0^{\infty}\hat{F}(ts)dF_T(t)\bigg)^mP_N\bigg(\int_0^{\infty}\hat{F}(ts)dF_T(t)\bigg)\\
&=&\big(\E[\exp(-sTX)]\big)^mP_N\big(\E[\exp(-sTX)]\big),\ \ s\ge 0,\nonumber
\end{eqnarray}
where $N\ge 0$ as in Eq (2).

In the next two theorems, we consider Eq (9) for the cases $m=1$ and
$m\ge 2,$ separately (see the explanations right after Theorem 3).

\noindent{\bf Theorem 2.} Let $0\le X\sim F$ with Laplace--Stieltjes
transform $\hat{F}$ and $\mu=\E[X]\in(0,\infty).$ Let $T\ge 0$ and
$N\ge 0$ be two given random variables, where $N$ takes values in
$\mathbb{N}_0$ and has pgf $P_N.$  Then for given $\mu,$ the random
variables $N$ and $T$ together satisfy the conditions
\begin{eqnarray} \E[N]=\frac{1-\E[T]}
{\E[T]},\ \ 0<\E[T^2]<\E[T]<1,\ \ \hbox{and}\ \ \E[N^2]<\infty
\end{eqnarray} iff the functional equation (9) with $m=1$  has
exactly one solution $F$ with mean $\mu$ and a finite variance.
Moreover, the variance is of the form
\begin{eqnarray} \hbox{Var}(X)=\frac{(\E[T])^2
\hbox{Var}(N)+\E[N+1] \hbox{Var}(T)}{1-\E[N+1]\E[T^2]} \cdot \mu^2
\end{eqnarray}
with $\E[N]=(1-\E[T])/{\E[T]}.$

 \noindent{\bf Corollary 3.} Under the setting of Theorem 2, the functional equation (9) with $m=1$ has
exactly one solution $F$ degenerate at mean $\mu=\E[X]$ iff the
random variables  $N$ and $T$ are degenerate at $\E[N]$ and $\E[T],$
respectively, and $\E[T]=1/(\E[N]+1)\in(0,1);$ precisely,
$\Pr(N=n_0)=1$ for some integer $n_0\ge 1$ and $\Pr(T=1/(n_0+1))=1.$

When $\Pr(T=p)=1$ for some $p\in(0,1),$ Theorem 2 reduces to the
following.

 \noindent{\bf Corollary 4.} Under the setting of Corollary 2, the random variable
 $N$ satisfies the conditions
 \[\E[N]=(1-p)/p\ \ \hbox{and}\ \ \E[N^2]<\infty\]
iff the functional equation
\begin{eqnarray*}\hat{F}(s)=\hat{F}(ps)P_N\big(\hat{F}(ps)\big),\ \ s\ge 0,\end{eqnarray*}
has exactly one solution $F$ with mean $\mu$ and a finite variance.
Moreover, the variance equals
\[ \hbox{Var}(X)=\frac{p^2 \hbox{Var}(N)}{1-p} \cdot \mu^2.\]

\noindent{\bf Theorem 3.} Let $0\le X\sim F$ with Laplace--Stieltjes
transform $\hat{F}$ and $\mu=\E[X]\in(0,\infty).$ Let $T\ge 0$ and
$N\ge 0$ be two given random variables, where $N$ takes values in
$\mathbb{N}_0$ and has pgf $P_N.$ Assume further that $m\ge 2$ is an
integer. Then for given $\mu,$ the random variables $N$ and $T$
together satisfy the conditions
\begin{eqnarray} \E[N]=\frac{1-m\E[T]}
{\E[T]},\ \ 0<\E[T^2]<\E[T]\le \frac1m<1,\ \ \hbox{and}\ \ \E[N^2]<\infty
\end{eqnarray} iff the functional equation (9) has exactly one
solution $F$ with mean $\mu$ and a finite variance. Moreover, the
variance is of the form
\begin{eqnarray} \hbox{Var}(X)=\frac{(\E[T])^2
\hbox{Var}(N)+\E[N+m] \hbox{Var}(T)}{1-\E[N+m]\E[T^2]} \cdot \mu^2
\end{eqnarray}
with $\E[N]=(1-m\E[T])/{\E[T]}.$
\medskip\\
\indent Note that we don't exclude the case $N=0$ in Theorem 3,
because when $N=0,$ Eq (9) with $m\ge 2$ is not a trivial case.
Besides, when $N=0,$ $\E[T]$ in (12) equals $1/m\,(\le 1/2<1),$
while the first two conditions in (10) fail to hold together.
Therefore, Theorem 2 is not a special case of Theorem 3; namely, we
cannot derive Theorem 2 from Theorem 3 by just letting $m=1.$ On the
other hand, it is seen that Eq (9) with $m\ge 2$ and $N=0$ is
equivalent to Eq (2) with $N=m\ge 2.$

 \noindent{\bf Corollary 5.} Under the setting of Theorem 3, the functional equation (9) has
exactly one solution $F$ degenerate at mean $\mu=\E[X]$ iff the
random variables  $N$ and $T$ are degenerate at $\E[N]$ and $\E[T],$
respectively, and $\E[T]=1/(\E[N]+m)\in(0,1/m];$ precisely,
$\Pr(N=n_0)=1$ for some integer $n_0\ge 0$ and $\Pr(T=1/(n_0+m))=1.$

When $\Pr(T=p)=1$ for some $p\in(0,1/m],$ Theorem 3 reduces to the
following.

 \noindent{\bf Corollary 6.} Under the setting of Corollary 2, assume, in addition,  $p\in(0,1/m],$
 where $m\ge 2$ is an
integer. Then for given $\mu,$ the random variable
 $N$ satisfies the conditions
 \[\E[N]=(1-mp)/p\ \ \hbox{and}\ \ \E[N^2]<\infty\]
iff the functional equation
\begin{eqnarray*}\hat{F}(s)=(\hat{F}(ps))^mP_N\big(\hat{F}(ps)\big),\ \ s\ge 0,\end{eqnarray*}
has exactly one solution $F$ with mean $\mu$ and a finite variance.
Moreover, the variance equals
\[ \hbox{Var}(X)=\frac{p^2 \hbox{Var}(N)}{1-p} \cdot \mu^2.\]

The distributional equations (1) and (8) are homogeneous cases,
because $X, X_1,X_2,\ldots$ are i.i.d random variables. We now
consider a nonhomogeneous case defined below. In addition to the
setting for Eq (1), suppose $0\le B\sim F_B$ is another random
variable independent of all $X, X_i, T, T_i, N.$ We will find
necessary and sufficient conditions  on $B,\ T$ and $N$ for which
the distributional equation
\begin{eqnarray}
X\stackrel{\rm d}{=}B+\sum_{i=1}^NT_iX_i,
\end{eqnarray}
 has a solution $X\sim F$ with
finite {\it variance}. Like Eq (2), Eq (14) has the  functional form
\begin{eqnarray}
\hat{F}(s)&=&\hat{F_B}(s) \cdot P_N\bigg(\int_0^{\infty}\hat{F}(ts)dF_T(t)\bigg)\\
&=&\hat{F_B}(s) \cdot P_N\big(\E[\exp(-sTX)]\big),\ \ \ s\ge 0.\nonumber
\end{eqnarray}

\noindent{\bf Theorem 4.} Let $0\le X\sim F$ with Laplace--Stieltjes
transform $\hat{F}$ and $\mu=\E[X]\in(0,\infty).$ Let $T\ge 0$ and
$N\ge 0$ be two given random variables with finite variances, where
$N$ takes values in $\mathbb{N}_0$ and has pgf $P_N.$ Suppose that
$0\le B\sim F_B$ is another random variable with mean $\E[B]>0$ and
a finite variance. Assume further that (i) $\Pr(N=0)<1,\ \Pr(T=0)<1$
and  (ii) $\hbox{Var}(B)+\hbox{Var}(T)+\hbox{Var}(N)>0.$ Then for
given $\mu,$ the following statements are true.
\\
(a) The random variables $B, N$ and $T$ together satisfy the
conditions
\begin{eqnarray} \mu=\frac{\E[B]}{1-\E[N]\E[T]},\ \ \ 0<\E[N]\E[T]<1,\ \ \ \hbox{and}\ \ \ 0<\E[N]\E[T^2]<1
\end{eqnarray} iff the functional equation
(15) has one solution $F$ with mean $\mu$ and a finite variance.
Moreover, the variance is of the form
\begin{eqnarray}
\hbox{Var}(X)=\frac{\hbox{Var}(B)+\mu^2(\E[T])^2
\hbox{Var}(N)+\mu^2\E[N]\hbox{Var}(T)}{1-\E[N]\E[T^2]}
\end{eqnarray}
with  $\mu={\E[B]}/(1-\E[N]\E[T]).$\\
(b) If, in addition to (16), $\E[T^2]<\E[T],$ then the solution $F$
to Eq (15) is unique.
\smallskip\\
\indent The purpose of the assumptions (i) and (ii) in Theorem 4  is
to exclude the trivial cases:\\
 (a) if $N=0$ or $T=0,$ Eq (14) reduces
to the equality $X\stackrel{\rm d}{=}B;$\\ (b) if
$\hbox{Var}(B)+\hbox{Var}(T)+\hbox{Var}(N)=0,$ all $B, T, N$ have
degenerate distributions, and so does the solution $X.$

The following interesting theorem points out the one-to-one
correspondence between solutions to Eq (1) and Eq (18) defined
below, where $\alpha\in(0,1).$

\noindent{\bf Theorem 5.} Let $\alpha\in(0,1)$ and let $0\le
T_{\alpha}\sim H_{\alpha}$ have the stable distribution with
Laplace--Stieltjes transform
$\hat{H}_{\alpha}(s)=\exp(-s^{\alpha}),\ s\ge 0.$ Then, under the
setting of Eq (1) with given $N$ and $T,$ \ $X_*\sim F_*$ is a
solution to Eq (1) with a mean $\mu\in(0,\infty)$ iff
$X_{\alpha}\sim F_{\alpha},$ where
$X_{\alpha}=T_{\alpha}X_*^{1/{\alpha}}$ and $T_{\alpha}$ is
independent of $X_*,$ is a solution to the distributional equation
\begin{eqnarray}
X\stackrel{\rm d}{=}\sum_{i=1}^NT_i^{1/{\alpha}}X_i
\end{eqnarray}
with $\lim_{s\to
0^+}(1-\hat{F}_{\alpha}(s))/s^{\alpha}=\mu\in(0,\infty).$
\bigskip\\
 \noindent{\bf 3. Lemmas}

 To prove the main results, we need some lemmas in the sequel.
Recall that  the pgf $P_N$ of a random variable $N$ taking values in
$\mathbb{N}_0$ is an absolutely monotone function on $[0,1]$ with
$P_N(1)=1,$ because $P_N(t)=\E[t^N]=\sum_{n=0}^{\infty}r_nt^n,\
t\in[0,1],$ with each $r_n=\Pr(N=n)\ge 0.$ For the first two lemmas,
see, e.g., Steutel and van Harn (2004), pp.\,483--484; Lemma 1 is
the so-called Bernstern Theorem.

\noindent{\bf Lemma 1.}  The Laplace--Stieltjes transform $\hat{F}$
of a nonnegative random variable $X\sim F$ is a completely monotone
function on $[0,\infty)$ with $\hat{F}(0)=1,$ and vice versa.

\noindent{\bf Lemma 2.} Let $Q$ be a pgf on $[0,1]$ and let
$\rho_1,\ \rho_2$
 be two completely monotone functions on $[0,\infty)$
with $\rho_1(0)=\rho_2(0)=1.$ Then  each of the  composition
function
 $Q\circ \rho_1$ and the product function $\rho_1 \rho_2$
  is completely
monotone on $[0,\infty),$ and is the  Laplace--Stieltjes transform
of a nonnegative random variable.

\noindent{\bf Lemma 3.} If $a, b\in[0,1]$ and $t\ge 1$ are three
real numbers, then $|a^t-b^t|\le t|a-b|.$\\
{\bf Proof.} If $t=1,$ the result is trivial. Suppose now that
$t>1.$ There are two possible cases for $a$ and $b$: (i) $0\le a \le
b\le 1$ and (ii) $0\le b \le a\le 1.$ It suffices to prove Case (i),
because Case (ii) follows from Case (i) immediately by the symmetry
property. Consider the function: $g(x)=x^t-tx,\ x\in[0,1].$ Since
$g^{\prime}(x)=t(x^{t-1}-1)\le 0,\ x\in[0,1],$ the function $g$ is
decreasing
 on $[0,1],$ Therefore, $g(a)\ge g(b)$ for Case (i). That
is, $a^t-ta\ge b^t-tb$ for $0\le a \le b\le 1.$ Equivalently,
$b^t-a^t\le t(b-a)$ \ or \[|a^t-b^t|=b^t-a^t\le t(b-a)= t|a-b|\ \
\hbox{for}\ \ 0\le a \le b\le 1.\] The proof is complete.

For a proof of the next crucial lemma, see Eckberg (1977), Guljas et
al.\,(1998) or Hu and Lin (2008).

\noindent{\bf Lemma 4.} Let $0\le X\sim F$ have a finite positive
second moment.  Then its Laplace--Stieltjes transform satisfies
\begin{eqnarray}
\hat{F}(s)\le
1-\frac{\mu_1^2}{\mu_2}+\frac{\mu_1^2}{\mu_2}e^{-(\mu_2/\mu_1)s},\ \
s\ge 0,
\end{eqnarray} where $\mu_j,\ j=1,2,$ is the $j$th moment of $X.$
\medskip\\
\indent
Note that in Lemma 4,  if the variance of $X$ is zero, then
$\mu_2=\mu_1^2$ and $X$ is degenerate at the mean $\mu_1>0.$ In this
case, (19) becomes an equality: $\hat{F}(s)=e^{-\mu_1s},\ s\ge 0.$

The next two lemmas are taken from Lin (2003, 1993). The sufficiency
parts of Corollaries 1, 3, and 5 can be proved directly by using
Lemma 5.

\noindent{\bf Lemma 5.} Let $g$ be a nonnegative function defined on
$[0,\infty)$ and let $g$ satisfy (i) $g(0)=1$, (ii)
$g^{\prime}(0)=b\in\mathbb{R}\equiv (-\infty,\infty)$ and (iii) for
some positive real $r\ne 1,$ $g(rx)=(g(x))^r,\ \ x\ge 0.$ Then $g$
is the exponential function $g(x)=e^{bx},\  x\ge 0.$

 \noindent{\bf Lemma 6.} Let
$0\le X\sim F$ with Laplace--Stieltjes transform $\hat{F}.$ Then for
each integer $n\ge 1,$ the $n$th moment $\E[X^n]=\lim_{s\to
0^+}(-1)^n\hat{F}^{(n)}(s)=(-1)^n\hat{F}^{(n)}(0^+)$ (finite or
infinite).

 For $0\le X\sim F$ with finite positive mean
$\mu_1,$ we define the first-order equilibrium distribution by
$F_{(1)}(x)=\mu_1^{-1}\int_0^x\overline{F}(y)\,dy,\ x\ge 0,$ where
$\overline{F}(x)=1-F(x).$ The high-order equilibrium distributions
are defined iteratively. Namely, the $n$th-order equilibrium
distribution  is
$F_{(n)}(x)=\mu_{(n-1)}^{-1}\int_0^x\overline{F}_{(n-1)}(y)\,dy,\
x\ge 0,$  provided the mean $\mu_{(n-1)}$ of ${F}_{(n-1)}$ is finite
(equivalently, the $n$th moment $\mu_n=\E[X^n]$ of $F$  is finite).
For the next relationship between the means of $\{F_{(n)}\}$ and
moments of $F$, see, e.g., Lin (1998), p.\,265, or Harkness and
Shantaram (1969).

\noindent{\bf Lemma 7.} Let  $0\le X\sim F$ have the $n$th moment
$\mu_n\in(0,\infty)$ for some $n\ge 2.$  Then the mean of the
$(n-1)$th-order equilibrium distribution $F_{(n-1)}$ is equal to
 $\mu_{(n-1)}=\mu_{n}/(n\,\mu_{n-1}).$

\noindent{\bf Lemma 8.} Let $0\le X\sim F$ with finite mean
$\mu\in(0,\infty)$ and let $X_{(1)}\sim F_{(1)}$ have the
first-order equilibrium  distribution.
 Then for $s>0,$ the following statements are true:\\
(i) $(1-\hat{F}(s))/s=\int_0^{\infty}e^{-sx}(1-F(x))dx;$\\
(ii) $\hat{F}_{(1)}(s)=(1-\hat{F}(s))/(\mu s)\le 1;$\\
(iii) $(\hat{F}(s)-1+\mu
s)/s^2=\mu\int_0^{\infty}e^{-sx}(1-F_{(1)}(x))dx;$\\
(iv) $\lim_{s\to 0^+}(1-\hat{F}(s))/s=\mu$ and $\lim_{s\to
0^+}(\hat{F}(s)-1+\mu s)/s^2=\E[X^2]/2$  (finite or infinite).\\
{\bf Proof.} For part (i), see Lin (1998), p.\,260, or Feller
(1971), p.\,435. Parts (ii) -- (iv) follow from the definition of
equilibrium distribution and Lemmas 6 and 7 immediately.

The next two lemmas are  key tools to prove the main results.

\noindent{\bf Lemma 9.} Let $0\le Y_n\sim G_n,$ $n=0,1,2,\ldots,$ be
a sequence of random variables having
 the same first two finite moments, say $\mu_1$ and
$\mu_2$. Suppose that their Laplace--Stieltjes transforms
$\{\hat{G}_n\}_{n=0}^{\infty}$ form a decreasing sequence of
functions. Then the limiting function
$\lim_{n\to\infty}\hat{G}_n(s)=\hat{G}_{\infty}(s),\ s\ge 0,$ exists
and is the Laplace--Stieltjes transform of a nonnegative random
variable, say $Y_{\infty},$ which has the mean
$\E[Y_{\infty}]=\mu_1$ and second
moment $\E[Y_{\infty}^2]\in[\mu_1^2,\mu_2].$\\
{\bf Proof.} For each fixed $s\ge 0,$ $\hat{G}_n(s)\in[0,1],\ n\ge
0,$ so the decreasing sequence $\{\hat{G}_n(s)\}_{n=0}^{\infty}$ has
a limit, denoted
$\hat{G}_{\infty}(s)=\lim_{n\to\infty}\hat{G}_n(s).$ On the other
hand, we have, by Jensen's inequality and the assumption,
\begin{eqnarray}e^{-\mu_1 s}\le
\hat{G}_{n}(s)\le \hat{G}_{0}(s),\ \ s\ge 0,\ n\ge 1.
\end{eqnarray}
 Therefore, the limiting function $\hat{G}_{\infty}$ satisfies
$e^{-\mu_1 s}\le \hat{G}_{\infty}(s)\le \hat{G}_{0}(s),\  s\ge 0,$
and hence $\lim_{s\to
0^+}\hat{G}_{\infty}(s)=1=\hat{G}_{\infty}(0).$ By the continuity
theorem for Laplace--Stieltjes transforms (see, e.g., Steutel and
van Harn 2004, p.\,479), we conclude that $\hat{G}_{\infty}$ is the
Laplace--Stieltjes transform of a nonnegative random variable,
denoted $Y_{\infty}.$ It remains to verify $\E[Y_{\infty}]=\mu_1$
and $\E[Y_{\infty}^2]\in[\mu_1^2,\mu_2].$  From (20) it follows that
the limiting function $\hat{G}_{\infty}$ satisfies
\begin{eqnarray}\frac{1-e^{-\mu_1 s}}{s}\ge
\frac{1-\hat{G}_{\infty}(s)}{s}\ge \frac{1-\hat{G}_{0}(s)}{s},\ \ s> 0,
\end{eqnarray}
and
\begin{eqnarray}\frac{e^{-\mu_1 s}-1+\mu_1 s}{s^2}\le
\frac{\hat{G}_{\infty}(s)-1+\mu_1 s}{s^2}\le \frac{\hat{G}_{0}(s)-1+\mu_1 s}{s^2},\ \ s> 0.
\end{eqnarray}
Finally, applying Lemma 8(iv) first to (21) gets
$\E[Y_{\infty}]=\mu_1$ and then to (22) yields $\mu_1^2/2\le
\E[Y_{\infty}^2]/2\le\mu_2/2.$ This completes the  proof.

 \noindent{\bf Lemma 10.} Let $W_1\sim F_{W_1}$ and
$W_2\sim F_{W_2}$ be two nonnegative random variables with the same
mean $\mu_W\in(0,\infty),$ and let $0\le Z^*\sim F_{Z^*}$ have a
mean $\mu_{Z^*}\in(0,1).$ Assume further that the Laplace--Stieltjes
transforms of $W_1$ and $W_2$ satisfy
\begin{eqnarray}
|\hat{F}_{W_1}(s)-\hat{F}_{W_2}(s)|\le\int_0^{\infty}|\hat{F}_{W_1}(ts)-\hat{F}_{W_2}(ts)|dF_{Z^*}(t),\
\ \ s\ge 0,
\end{eqnarray} or, equivalently,
\[\big|\E[\exp(-sW_1)]-\E[\exp(-sW_2)]\big|\le
\big|\E[\exp(-sZ^*W_1)]-\E[\exp(-sZ^*W_2)]\big|,\ \ s\ge 0.\] Then
$\hat{F}_{W_1}=\hat{F}_{W_2}$ and hence $F_{W_1}=F_{W_2}.$\\
{\bf Proof.} For  each $s>0,$ applying the inequality (23) $(n-1)$
more times yields
\begin{eqnarray}
& &|\hat{F}_{W_1}(s)-\hat{F}_{W_2}(s)|\le\int_0^{\infty}\!\!\cdots\!\!\int_0^{\infty}|\hat{F}_{W_1}(t_1\cdots t_ns)-
\hat{F}_{W_2}(t_1\cdots t_ns)|dF_{Z^*}(t_1)\cdots dF_{Z^*}(t_n)\nonumber\\
&= &\int_{0^+}^{\infty}\!\!\cdots\!\!\int_{0^+}^{\infty}\bigg|\frac{\hat{F}_{W_1}(t_1\cdots t_ns)-\hat{F}_{W_2}(t_1\cdots t_ns)}{
\mu_Wt_1\cdots t_ns}\bigg|\mu_Wt_1\cdots t_nsdF_{Z^*}(t_1)\cdots dF_{Z^*}(t_n).
\end{eqnarray}
We now estimate the first part of the integrand:
\begin{eqnarray}
& &\bigg|\frac{\hat{F}_{W_1}(t_1\cdots t_ns)-\hat{F}_{W_2}(t_1\cdots t_ns)}{
\mu_Wt_1\cdots t_ns}\bigg|=\bigg|\frac{1-\hat{F}_{W_1}(t_1\cdots t_ns)}{\mu_Wt_1\cdots t_ns}-\frac{1-\hat{F}_{W_2}(t_1\cdots t_ns)}{
\mu_Wt_1\cdots t_ns}\bigg|\nonumber\\
& \le &\bigg|\frac{1-\hat{F}_{W_1}(t_1\cdots t_ns)}{\mu_Wt_1\cdots t_ns}\bigg|+\bigg|\frac{1-\hat{F}_{W_2}(t_1\cdots t_ns)}{
\mu_Wt_1\cdots t_ns}\bigg|\,\ \le\ 2.
\end{eqnarray}
The last inequality is due to Lemma 8(ii). Combining (24) and (25)
together  leads to
\begin{eqnarray*}|\hat{F}_{W_1}(s)-\hat{F}_{W_2}(s)|&\le &
2\mu_Ws\int_{0^+}^{\infty}\!\!\cdots\!\!\int_{0^+}^{\infty}t_1\cdots
t_ndF_{Z^*}(t_1)\cdots dF_{Z^*}(t_n)\\
&=&2\mu_Ws(\E[Z^*])^n\ \ \longrightarrow 0\ \ \ \ \hbox{as}\ \ n\to \infty,
\end{eqnarray*}
in which the last conclusion follows from the assumption
$\E[Z^*]=\mu_{Z^*}\in (0,1).$ Therefore,
$\hat{F}_{W_1}=\hat{F}_{W_2}.$ This completes the proof.
\bigskip\\
 \noindent{\bf 4. Proofs of main results}

\noindent{\bf Proof of Theorem 1.} (Sufficiency) Suppose that Eq (2)
has exactly one solution $0\le X\sim F$ with mean $\mu\in(0,\infty)$
and a finite variance (and hence $\E[X^2]\in(0,\infty)$). Then we
want to prove that the conditions (6) hold true.

Rewrite Eq (2) as
\[\hat{F}(s)=P_N\bigg(\int_0^{\infty}\hat{F}(ts)dF_T(t)\bigg)=\sum_{n=0}^{\infty}\Pr(N=n)\bigg(\int_0^{\infty}\hat{F}(ts)dF_T(t)\bigg)^n,\ \ s \ge 0.\]
Differentiating twice the above equation with respect to $s,$ we
have, for $s>0,$
\begin{eqnarray}
\hat{F}^{\prime}(s)&=&\sum_{n=1}^{\infty}\Pr(N=n)n\bigg(\int_0^{\infty}\hat{F}(ts)dF_T(t)\bigg)^{n-1}\int_0^{\infty}\hat{F}^{\prime}(ts)tdF_T(t),\\
\hat{F}^{\prime\prime}(s)&=&\sum_{n=2}^{\infty}\Pr(N=n)n(n-1)\bigg(\int_0^{\infty}\hat{F}(ts)dF_T(t)\bigg)^{n-2}\bigg(\int_0^{\infty}\hat{F}^{\prime}(ts)tdF_T(t)\bigg)^2\nonumber\\
& &+\sum_{n=1}^{\infty}\Pr(N=n)n\bigg(\int_0^{\infty}\hat{F}(ts)dF_T(t)\bigg)^{n-1}\int_0^{\infty}\hat{F}^{\prime\prime}(ts)t^2dF_T(t).
\end{eqnarray}
Letting $s\to 0^+$ in (26) and (27) yields, respectively,
\begin{eqnarray*}
\hat{F}^{\prime}(0^+)&=&\hat{F}^{\prime}(0^+)\E[N]\E[T],\\
\hat{F}^{\prime\prime}(0^+)&=&\E[N(N-1)](\hat{F}^{\prime}(0^+)\E[T])^2+\hat{F}^{\prime\prime}(0^+)\E[N]\E[T^2].
\end{eqnarray*}
Equivalently, we have, by Lemma 6,
\begin{eqnarray}
\mu&=&\mu\ \E[N]\E[T],\\
\E[X^2]&=&\E[N(N-1)](\mu\,\E[T])^2+\E[X^2]\E[N]\E[T^2].
\end{eqnarray}
From (28) and (29) it follows that $\E[N]\E[T]=1$ (which implies
that $\E[N], \E[T]>0$) and $\E[N^2]<\infty$ because $\mu,
\E[X^2]\in(0,\infty).$ It remains to prove that $0<\E[T^2]<\E[T]<1.$

Since $\E[N(N-1)]\ge 0,$ we have by (29)  that $\E[N]\E[T^2]\le 1,$
and hence $\E[T^2]\le\E[T]$ due to the fact $\E[N]\E[T]=1.$ Namely,
\begin{eqnarray}0<(\E[T])^2\le \E[T^2]\le \E[T],
\end{eqnarray}
from which we further have $0<\E[T]\le 1.$ We now prove that
$\E[T]<1.$ Suppose on the contrary $\E[T]=1.$ Then $\E[N]=1$ (by the
fact $\E[N]\E[T]=1$) and from (30) it follows that $\E[T^2]=1,\
\hbox{Var}(T)=0$ and $\Pr(T=1)=1.$ Plugging these in (29) yields
\begin{eqnarray}1=(\E[N])^2\le \E[N^2]= \E[N]=1,
\end{eqnarray}
which implies $\Pr(N=1)=1$ as in the case of $T.$ These together
imply that Eq (2) is an identity for any $0\le X\sim F$ as described
before, which contradicts the unique solution to Eq (2). So we
conclude that $\E[T]\in(0,1).$

Finally, we prove  $\E[T^2]<\E[T].$  Suppose on the contrary
$\E[T^2]=\E[T].$ Then (31) follows from (29) again (using
$\E[N]\E[T]=1$) and hence $\Pr(N=1)=1$. This is impossible because
$\E[T]\in(0,1)$ and $\E[N]\E[T]=1.$ The proof of the sufficiency
part is complete.

 (Necessity) Suppose that the conditions (6) hold true. Then we will prove the existence of
 a solution $F$ to Eq (2) with mean $\mu$ and a finite variance.

 Set first
 \begin{eqnarray}\mu_1=\mu\ \ \hbox{and}\ \
 \mu_2=\frac{\E[N(N-1)](\E[T])^2}{1-\E[N]\E[T^2]}\cdot \mu_1^2.\end{eqnarray}
Note that the denominator $1-\E[N]\E[T^2]=1-\E[T^2]/\E[T]>0$ by (6)
and that $\mu_2\ge\mu_1^2$ by the facts: $\E[N^2]\ge(\E[N])^2$ and
$\E[T^2]\ge(\E[T])^2.$ Therefore, the RHS of (19) with $\mu_1,
\mu_2$ defined in  (32) is a bona fide Laplace--Stieltjes transform,
say $\hat{F}_0,$ of a nonnegative random variable $Y_0\sim F_0$ (by
Lemma 1). Namely,
\[\hat{F}_0(s)=1-\frac{\mu_1^2}{\mu_2}+\frac{\mu_1^2}{\mu_2}e^{-(\mu_2/\mu_1)s},\ \
s\ge 0.
\]

Next, using the initial $Y_0\sim F_0$ we define iteratively the
sequence of random variables $Y_n\sim F_n,\ n=1,2,\ldots,$ through
Laplace--Stieltjes transforms:
\begin{eqnarray}
\hat{F}_n(s)=P_N\bigg(\int_0^{\infty}\hat{F}_{n-1}(ts)dF_T(t)\bigg)=P_N\big(\E[\exp(-sTY_{n-1})]\big),\ n\ge 1,
\end{eqnarray}
which is well-defined due to Lemma 2. Differentiating twice the
above equation with respect to $s$ and letting $s\to 0^+,$ we have,
for $n\ge 1,$
 \begin{eqnarray}
\hat{F}_n^{\prime}(0^+)&=&\hat{F}_{n-1}^{\prime}(0^+)\E[N]\E[T]=\hat{F}_{n-1}^{\prime}(0^+),\\
 \hat{F}_n^{\prime\prime}(0^+)&=&\E[N(N-1)](\hat{F}_{n-1}^{\prime}(0^+)\E[T])^2+\hat{F}_{n-1}^{\prime\prime}(0^+)\E[N]\E[T^2].\end{eqnarray}

With the help of Lemma 6 and  by induction on $n,$ we can show
through (34) and (35) that $\E[Y_n]=\E[Y_0]=\mu_1,$
$\E[Y_n^2]=\E[Y_0^2]=\mu_2$ (defined in (32)) for all $n\ge 1$ and
hence
\begin{eqnarray}\hbox{Var}(Y_n)=\mu_2-\mu_1^2=\frac{(\E[T])^2\hbox{Var}(N)+\E[N]\hbox{Var}(T)}{1-\E[N]\E[T^2]}\cdot
\mu_1^2,\ \ n\ge 0.
\end{eqnarray}
Moreover, by Lemma 4, we first have $\hat{F}_1\le \hat{F}_0$, and
then by the iteration (33), $\hat{F}_n\le \hat{F}_{n-1}$ for all
$n\ge 2$ (due to the absolute monotonicity of $P_N$). Namely,
$\{Y_n\}_{n=0}^{\infty}$ is a sequence of nonnegative random
variables having the same first two moments $\mu_1, \mu_2,$ and
their Laplace--Stieltjes transforms $\{\hat{F}_n\}$ are decreasing.
Therefore, Lemma 9 applies. Denote the limit of $\{\hat{F}_n\}$ by
$\hat{F}_{\infty},$ which is the Laplace--Stieltjes transform of a
nonnegative random variable $Y_{\infty}\sim F_{\infty}$ with
$\E[Y_{\infty}]=\mu_1$ and $\E[Y_{\infty}^2]\in[\mu_1^2,\mu_2].$
Consequently, it follows from (33) that the limit $F_{\infty}$ is a
solution to Eq (2) with mean
 $\mu$ and a
finite variance. Applying Lemma 6 to Eq (2) again (with
$F=F_{\infty}$), we conclude that $\E[Y_{\infty}^2]=\mu_2$ as given
in (32), and hence the solution ${Y}_{\infty}\sim F_{\infty}$ has
the required variance as shown in (7) or (36).

Finally, we prove the uniqueness of the solution to Eq (2). Suppose
 there are two solutions to Eq (2), denoted $0\le X\sim F_X$ and $0\le Y\sim
F_Y.$  We want to show that $F_X=F_Y.$ As before, with the help of
Lemma 6 and the conditions (6), we have from Eq (2) that
\[\E[X]=\E[Y]=\mu_1=\mu,\ \ \ \E[X^2]=\E[Y^2]=\mu_2=\frac{\E[N(N-1)](\E[T])^2}{1-\E[N]\E[T^2]}\cdot \mu_1^2.\]
Let $W_1\sim F_{W_1}, W_2\sim F_{W_2}$ have the first-order
equilibrium distributions of $X, Y,$ respectively. By Lemma 8(ii),
their Laplace-Stieltjes transforms are of the form:
\begin{eqnarray}\hat{F}_{W_1}(s)=\frac{1-\hat{F}_X(s)}{\mu s},\quad
\hat{F}_{W_2}(s)=\frac{1-\hat{F}_Y(s)}{\mu s},\ s>0.
\end{eqnarray}
Therefore, it remains to prove that
$\hat{F}_{W_1}(s)=\hat{F}_{W_2}(s),\ s>0.$

From Lemma 7 it follows that
$\E[W_1]=\E[W_2]={\mu_2}/{(2\mu_1)}\equiv \mu_W\in(0,\infty).$ Using
Eq (2), we first estimate the difference between $\hat{F}_{X}$ and
$\hat{F}_{Y}$ as follows: for $s>0,$
\begin{eqnarray*}
& &|\hat{F}_{X}(s)-\hat{F}_{Y}(s)|=\bigg| P_N\bigg(\int_0^{\infty}\hat{F_X}(ts)dF_T(t)\bigg)-P_N\bigg(\int_0^{\infty}\hat{F_Y}(ts)dF_T(t)\bigg)\bigg|\\
&=&\bigg|\sum_{n=0}^{\infty}\Pr(N=n)\bigg[\bigg(\int_0^{\infty}\hat{F_X}(ts)dF_T(t)\bigg)^n-\bigg(\int_0^{\infty}\hat{F_Y}(ts)dF_T(t)\bigg)^n\bigg]\bigg|\\
&\le &\sum_{n=0}^{\infty}\Pr(N=n)\cdot n\int_0^{\infty}|\hat{F}_{X}(ts)-\hat{F}_{Y}(ts)|dF_T(t)\\
&=&\E[N]\int_0^{\infty}|\hat{F}_{X}(ts)-\hat{F}_{Y}(ts)|dF_T(t)\\
&=&\frac{1}{\E[T]}\int_0^{\infty}|\hat{F}_{X}(ts)-\hat{F}_{Y}(ts)|dF_T(t),
\end{eqnarray*}
in which the inequality follows from Lemma 3, while the last
equality is due to the condition $\E[N]\E[T]=1$ in (6). Therefore,
we have, for $s>0,$
\begin{eqnarray*}
\bigg|\frac{\hat{F}_{X}(s)-\hat{F}_{Y}(s)}{\mu s}\bigg|\le \int_0^{\infty}\bigg|\frac{\hat{F}_{X}(ts)-\hat{F}_{Y}(ts)}{\mu st}\bigg|\frac{t}{\E[T]}dF_T(t)
\equiv \int_0^{\infty}\bigg|\frac{\hat{F}_{X}(ts)-\hat{F}_{Y}(ts)}{\mu st}\bigg|dF_{Z^*}(t),
\end{eqnarray*}
where $Z^*\sim F_{Z^*}$ has the length-biased distribution of $T\sim
F_T$ and $\E[Z^*]=\E[T^2]/\E[T]<1.$ Equivalently, we have,  by (37),
that
\[|\hat{F}_{W_1}(s)-\hat{F}_{W_2}(s)|\le \int_0^{\infty}|\hat{F}_{W_1}(ts)-\hat{F}_{W_2}(ts)|dF_{Z^*}(t),\ s>0.
\]
 Lemma 10 applies and hence
$\hat{F}_{W_1}=\hat{F}_{W_2}.$ This proves the uniqueness of the
solution to Eq (2). The proof of the necessity part is complete.

\noindent{\bf Proof of Proposition 1.} (i) Since $\E[N]\E[T]=1,$ we
have $\E[N],\ \E[T]>0$ and hence $\Pr(T>0)>0.$ Rewrite
\[1=\E[N]\E[T]=\E[N]\E[T|\,T>0]\Pr(T>0).
\]
Therefore, $0<\Pr(T>0)\E[N]=(\E[T|\, T>0])^{-1}.$ We will show that
$\E[T|\, T>0])<1.$ Write
\[\E[T|\, T>0]=\frac{1}{\Pr(T>0)}\int_{0^+}^{\infty}tdF_T(t)=\frac{\E[T]}{\Pr(T>0)}.
\]
Similarly, $\E[T^2|\, T>0]={\E[T^2]}/{\Pr(T>0)}.$ From the condition
$\E[T^2]<\E[T]$ it then follows that $(\E[T|\, T>0])^2\le \E[T^2|\,
T>0]<\E[T|\, T>0].$ Consequently, $\E[T|\, T>0]<1.$ This proves the
first conclusion of the proposition.\\
(ii) We next prove the second conclusion $\E[T\log T]<0.$ Note first
that the function $g(t)=t^2-t-t\log t\ge 0$ \ for $t>0.$ So we have
$\E[T^2-T-T\log T\,|\,T>0]\ge 0.$ Equivalently,
\[\E[T^2-T\,|\,T>0]\ge \E[T\log T\,|\,T>0].
\]
Finally, by the condition $0<\E[T^2]<\E[T]<1,$ we have
\begin{eqnarray*}
0>\E[T^2-T]&=&\E[T^2-T\,|\,T>0]\Pr(T>0)\\
&\ge & \E[T\log
T\,|\,T>0]\Pr(T>0)=\E[T\log T].
\end{eqnarray*}
This completes the proof.

\noindent{\bf Proof of Theorem 2.}  Note that Eq (9) with $m=1$ is
equivalent to
\[\hat{F}(s)=P_{N+1}\bigg(\int_0^{\infty}\hat{F}(ts)dF_T(t)\bigg)=P_{N+1}\big(\E[\exp(-sTX)]\big),\ \ s\ge 0,
\]
and that $\hbox{Var}(N+1)=\hbox{Var}(N)$ in (11). Therefore, Theorem
2 follows from Theorem 1 by replacing  $N$ by $N+1$ taking values in
$\mathbb {N}\equiv \{1,2,3,\ldots\}.$ The proof is complete.

\noindent{\bf Proof of Theorem 3.} The proof is similar to that of
Theorem 1. We give the details here for completeness.

(Sufficiency) Suppose that Eq (9) with $m\ge 2$ has exactly one
solution $0\le X\sim F$ with mean $\mu\in(0,\infty)$ and a finite
variance (and hence $\E[X^2]\in(0,\infty)$). Then we want to prove
that the conditions (12) hold true.

Rewrite Eq (9) with $m\ge 2$ as
\begin{eqnarray*}\hat{F}(s)&=&\bigg(\int_0^{\infty}\hat{F}(ts)dF_T(t)\bigg)^mP_N\bigg(\int_0^{\infty}\hat{F}(ts)dF_T(t)\bigg)\\
&=&\sum_{n=0}^{\infty}\Pr(N=n)\bigg(\int_0^{\infty}\hat{F}(ts)dF_T(t)\bigg)^{n+m},\ \ \ s\ge 0.\end{eqnarray*}
Differentiating twice the above equation with respect to $s,$ we
have, for $s>0,$
\begin{eqnarray}
\hat{F}^{\prime}(s)&=&\sum_{n=0}^{\infty}\Pr(N=n)(n+m)\bigg(\int_0^{\infty}\hat{F}(ts)dF_T(t)\bigg)^{n+m-1}\int_0^{\infty}\hat{F}^{\prime}(ts)tdF_T(t),\\
\hat{F}^{\prime\prime}(s)&=&\sum_{n=0}^{\infty}\Pr(N=n)(n+m)(n+m-1)
\bigg(\int_0^{\infty}\hat{F}(ts)dF_T(t)\bigg)^{n+m-2}\bigg(\int_0^{\infty}\hat{F}^{\prime}(ts)tdF_T(t)\bigg)^2\nonumber\\
& &+\sum_{n=0}^{\infty}\Pr(N=n)(n+m)\bigg(\int_0^{\infty}\hat{F}(ts)dF_T(t)\bigg)^{n+m-1}\int_0^{\infty}\hat{F}^{\prime\prime}(ts)t^2dF_T(t).
\end{eqnarray}
Letting $s\to 0^+$ in (38) and (39) yields, respectively,
\begin{eqnarray*}
\hat{F}^{\prime}(0^+)&=&\hat{F}^{\prime}(0^+)\E[N+m]\E[T],\\
\hat{F}^{\prime\prime}(0^+)&=&\E[(N+m)(N+m-1)](\hat{F}^{\prime}(0^+)\E[T])^2+\hat{F}^{\prime\prime}(0^+)\E[N+m]\E[T^2].
\end{eqnarray*}
Equivalently, we have, by Lemma 6,
\begin{eqnarray}
\mu&=&\mu\ \E[N+m]\E[T],\\
\E[X^2]&=&\E[(N+m)(N+m-1)](\mu\,\E[T])^2+\E[X^2]\E[N+m]\E[T^2].
\end{eqnarray}
Therefore, from (40) and (41) it follows that $\E[N+m]\E[T]=1$ and
$\E[N^2]<\infty$ because $\mu, \E[X^2]\in(0,\infty).$ It remains to
prove that $0<\E[T^2]<\E[T]\le 1/m<1.$

Since $\E[(N+m)N]\ge 0,$  it further follows from (41) that
$\E[N+m]\E[T^2]\le 1,$ and hence $\E[T^2]\le\E[T]$ because
$\E[N+m]\E[T]=1.$ The latter also implies that $\E[T]>0,$ and hence
\begin{eqnarray*}0<(\E[T])^2\le \E[T^2]\le \E[T]=\frac{1}{\E[N+m]}\le\frac{1}{m}<1.
\end{eqnarray*}

Finally, we prove  $\E[T^2]<\E[T].$  Suppose on the contrary
$\E[T^2]=\E[T].$ Then from (41) it follows (by using
$\E[N+m]\E[T]=1$) that
\[\E[(N+m)(N+m-1)]=0.\]
This is impossible because $N\ge 0$ and $m\ge 2.$  The proof of the
sufficiency part is complete.

 (Necessity) Suppose that the conditions (12) hold true. Then we will prove the existence of
 a solution $F$ to Eq (9)  with $m\ge 2,$ which has mean $\mu$ and a finite variance.

 Set first
 \begin{eqnarray}\mu_1=\mu\ \ \hbox{and}\ \
 \mu_2=\frac{\E[(N+m)(N+m-1)](\E[T])^2}{1-\E[N+m]\E[T^2]}\cdot \mu_1^2.\end{eqnarray}
Note that the denominator $1-\E[N+m]\E[T^2]=1-\E[T^2]/\E[T]>0$ by
(12) and that $\mu_2\ge\mu_1^2$ by the facts:
$\E[(N+m)^2]\ge(\E[N+m])^2$ and $\E[T^2]\ge(\E[T])^2.$ Therefore,
the RHS of (19) with $\mu_1, \mu_2$ defined in (42) is a bona fide
Laplace--Stieltjes transform, say $\hat{F}_0,$ of a nonnegative
random variable $Y_0\sim F_0.$ Namely,
\[\hat{F}_0(s)=1-\frac{\mu_1^2}{\mu_2}+\frac{\mu_1^2}{\mu_2}e^{-(\mu_2/\mu_1)s},\ \ \
s\ge 0.
\]

Next, using the initial $Y_0\sim F_0$ we define iteratively the
sequence of random variables $Y_n\sim F_n,\ n=1,2,\ldots,$ through
Laplace--Stieltjes transforms:
\begin{eqnarray}
\hat{F}_n(s)=\bigg(\int_0^{\infty}\hat{F}_{n-1}(ts)dF_T(t)\bigg)^mP_N\bigg(\int_0^{\infty}\hat{F}_{n-1}(ts)dF_T(t)\bigg),\ \ n\ge 1,
\end{eqnarray}
which is well-defined due to Lemma 2. Differentiating twice the
above equation with respect to $s$ and letting $s\to 0^+,$ we have,
for $n\ge 1,$
 \begin{eqnarray}
\hat{F}_n^{\prime}(0^+)&=&\hat{F}_{n-1}^{\prime}(0^+)\E[N+m]\E[T]=\hat{F}_{n-1}^{\prime}(0^+),\\
 \hat{F}_n^{\prime\prime}(0^+)&=&\E[(N+m)(N+m-1)](\hat{F}_{n-1}^{\prime}(0^+)\E[T])^2+\hat{F}_{n-1}^{\prime\prime}(0^+)\E[N+m]\E[T^2].\end{eqnarray}

With the help of Lemma 6 and  by induction on $n,$ we can show
through (44) and (45) that $\E[Y_n]=\E[Y_0]=\mu_1,$
$\E[Y_n^2]=\E[Y_0^2]=\mu_2$ (defined in (42)) for all $n\ge 1$ and
hence
\begin{eqnarray}\hbox{Var}(Y_n)=\mu_2-\mu_1^2=\frac{(\E[T])^2\hbox{Var}(N)+\E[N+m]\hbox{Var}(T)}{1-\E[N+m]\E[T^2]}\cdot
\mu_1^2,\ \ \ n\ge 0.
\end{eqnarray}
Moreover, by Lemma 4, we first have $\hat{F}_1\le \hat{F}_0$, and
then by the iteration (43), $\hat{F}_n\le \hat{F}_{n-1}$ for all
$n\ge 2.$ Namely, $\{Y_n\}_{n=0}^{\infty}$ is a sequence of
nonnegative random variables having the same first two moments
$\mu_1, \mu_2,$ and their Laplace--Stieltjes transforms
$\{\hat{F}_n\}$ are decreasing. Therefore, Lemma 9 applies. Denote
the limit of $\{\hat{F}_n\}$ by $\hat{F}_{\infty},$ which is the
Laplace--Stieltjes transform of a nonnegative random variable
$Y_{\infty}\sim F_{\infty}$ with $\E[Y_{\infty}]=\mu_1$ and
$\E[Y_{\infty}^2]\in[\mu_1^2,\mu_2].$ Consequently, it follows from
(43) that the limit $F_{\infty}$ is a solution to Eq (9) with $m\ge
2,$ which has mean
 $\mu$ and a
finite variance. Applying Lemma 6 to Eq (9) with $m\ge 2$ again, we
conclude that $\E[Y_{\infty}^2]=\mu_2$ as given in (42), and hence
the solution ${Y}_{\infty}\sim F_{\infty}$ has the required variance
as shown in (13) or (46).

Finally, we prove the uniqueness of the solution to Eq (9) with
$m\ge 2.$ Suppose that there are two solutions, denoted $0\le X\sim
F_X$ and $0\le Y\sim F_Y.$ We want to show  $F_X=F_Y.$ As before,
with the help of Lemma 6 and the conditions (12), we have from Eq
(9) with $m\ge 2$ that
\[\E[X]=\E[Y]=\mu_1=\mu,\ \ \ \E[X^2]=\E[Y^2]=\mu_2=\frac{\E[(N+m)(N+m-1)](\E[T])^2}{1-\E[N+m]\E[T^2]}\cdot \mu_1^2.\]
Let $W_1\sim F_{W_1}, W_2\sim F_{W_2}$ have the first-order
equilibrium distributions of $X, Y,$ respectively. By Lemma 8(ii),
their Laplace--Stieltjes transforms are of the form:
\begin{eqnarray}\hat{F}_{W_1}(s)=\frac{1-\hat{F}_X(s)}{\mu s},\quad
\hat{F}_{W_2}(s)=\frac{1-\hat{F}_Y(s)}{\mu s},\ \ s>0.
\end{eqnarray}
Therefore, it remains to prove that
$\hat{F}_{W_1}(s)=\hat{F}_{W_2}(s),\ s>0.$

From Lemma 7 it follows that
$\E[W_1]=\E[W_2]={\mu_2}/{(2\mu_1)}\equiv \mu_W\in(0,\infty).$ Using
Eq (9) with $m\ge 2$, we first estimate the difference between
$\hat{F}_{X}$ and $\hat{F}_{Y}$ as follows: for $s>0,$
\begin{eqnarray*}
|\hat{F}_{X}(s)-\hat{F}_{Y}(s)|
&=&\bigg|\sum_{n=0}^{\infty}\Pr(N=n)\bigg[\bigg(\int_0^{\infty}\hat{F_X}(ts)dF_T(t)\bigg)^{n+m}-\bigg(\int_0^{\infty}\hat{F_Y}(ts)dF_T(t)\bigg)^{n+m}\bigg]\bigg|\\
&\le &\sum_{n=0}^{\infty}\Pr(N=n)\cdot (n+m)\int_0^{\infty}|\hat{F}_{X}(ts)-\hat{F}_{Y}(ts)|dF_T(t)\\
&=&\E[N+m]\int_0^{\infty}|\hat{F}_{X}(ts)-\hat{F}_{Y}(ts)|dF_T(t)\\
&=&\frac{1}{\E[T]}\int_0^{\infty}|\hat{F}_{X}(ts)-\hat{F}_{Y}(ts)|dF_T(t),
\end{eqnarray*}
in which the inequality follows from Lemma 3, while the last
equality is due to the condition $\E[N+m]\E[T]=1$ in (12).
Therefore, we have, for $s>0,$
\begin{eqnarray*}
\bigg|\frac{\hat{F}_{X}(s)-\hat{F}_{Y}(s)}{\mu s}\bigg|\le \int_0^{\infty}\bigg|\frac{\hat{F}_{X}(ts)-\hat{F}_{Y}(ts)}{\mu st}\bigg|\frac{t}{\E[T]}dF_T(t)
\equiv \int_0^{\infty}\bigg|\frac{\hat{F}_{X}(ts)-\hat{F}_{Y}(ts)}{\mu st}\bigg|dF_{Z^*}(t),
\end{eqnarray*}
where $Z^*\sim F_{Z^*}$ has the length-biased distribution of $T\sim
F_T$ and $\E[Z^*]=\E[T^2]/\E[T]<1.$ Equivalently, we have,  by (47),
that
\[|\hat{F}_{W_1}(s)-\hat{F}_{W_2}(s)|\le \int_0^{\infty}|\hat{F}_{W_1}(ts)-\hat{F}_{W_2}(ts)|dF_{Z^*}(t),\ \ s>0.
\]
 Lemma 10 applies and hence
$\hat{F}_{W_1}=\hat{F}_{W_2}.$ This proves the uniqueness of the
solution to Eq (9) with $m\ge 2.$ The proof of the necessity part is
complete.

\noindent{\bf Proof of Theorem 4.} We first sketch the proof of part
(a).

(Sufficiency) Suppose that Eq (15)  has one solution $0\le X\sim F$
with mean $\mu\in(0,\infty)$ and a finite variance (and hence
$\E[X^2]\in(0,\infty)$). Then we want to prove that the conditions
(16) hold true.

Rewrite Eq (15)  as
\begin{eqnarray*}\hat{F}(s)=\hat{F}_B(s)\sum_{n=0}^{\infty}\Pr(N=n)\bigg(\int_0^{\infty}\hat{F}(ts)dF_T(t)\bigg)^{n},\ \ \ s \ge 0.\end{eqnarray*}
Differentiating twice the above equation with respect to $s$ and
letting $s\to 0^+$ yield
\begin{eqnarray*}
\hat{F}^{\prime}(0^+)&=&\hat{F}_B^{\prime}(0^+)+\hat{F}^{\prime}(0^+)\E[N]\E[T],\\
\hat{F}^{\prime\prime}(0^+)&=&\hat{F}_B^{\prime\prime}(0^+)+2\hat{F}_B^{\prime}(0^+)\hat{F}^{\prime}(0^+)\E[N]\E[T]\\
& & +\, \E[N(N-1)](\hat{F}^{\prime}(0^+)\E[T])^2+\hat{F}^{\prime\prime}(0^+)\E[N]\E[T^2].
\end{eqnarray*}
Equivalently, we have, by Lemma 6,
\begin{eqnarray}
\mu&=&\E[B]+\mu\ \E[N]\E[T]>\mu\ \E[N]\E[T],\\
\E[X^2]&=&\E[B^2]+2\mu\,\E[B]\E[N]\E[T]+\E[N(N-1)](\mu\,\E[T])^2+\E[X^2]\E[N]\E[T^2]\nonumber\\
&=&\hbox{Var}(B)+\mu^2(\E[T])^2\hbox{Var}(N)+\mu^2(1-\E[N](\E[T])^2)+\E[X^2]\E[N]\E[T^2].
\end{eqnarray}
Therefore, $0< \E[N]\E[T]<1$ by (48) and $\mu=\E[B]/(1-\E[N]\E[T]).$
Moreover, it follows from (48) and (49) that
\begin{eqnarray*}
\hbox{Var}(X)(1-\E[N]\E[T^2])=\hbox{Var}(B)+\mu^2(\E[T])^2\hbox{Var}(N)+\mu^2\E[N]\hbox{Var}(T)\in(0,\infty).
\end{eqnarray*}
The last conclusion is due to the assumptions. Thus, $0<
\E[N]\E[T^2]<1,$ and we have the variance of $X$ as required in
(17). The proof of the sufficiency part of (a) is complete.

 (Necessity) Suppose that the conditions (16) hold true. Then we will prove the existence of
 a solution $F$ to Eq (15),  which has mean $\mu$ and a finite variance.

 Set first
 \begin{eqnarray}\mu_1&=&\mu=\frac{\E[B]}{1-\E[N]\E[T]},\\
 \mu_2&=&\frac{\hbox{Var}(B)+\mu^2(\E[T])^2\hbox{Var}(N)+\mu^2(1-\E[N](\E[T])^2)}{1-\E[N]\E[T^2]}.\end{eqnarray}
Note that the denominators $1-\E[N]\E[T],\ 1-\E[N]\E[T^2]>0$ by (16)
and that $\mu_2\ge \mu_1^2$ by the fact $\hbox{Var}(T)\ge 0.$
Therefore,  the RHS of (19) with (50) and (51) is a bona fide
Laplace--Stieltjes transform, say $\hat{F}_0,$ of a nonnegative
random variable $Y_0\sim F_0.$ Namely,
\[\hat{F}_0(s)=1-\frac{\mu_1^2}{\mu_2}+\frac{\mu_1^2}{\mu_2}e^{-(\mu_2/\mu_1)s},\ \ \
s\ge 0.
\]

Next, using the initial $Y_0\sim F_0$ we define iteratively the
sequence of random variables $Y_n\sim F_n,\ n=1,2,\ldots,$ through
Laplace--Stieltjes transforms:
\begin{eqnarray}
\hat{F}_n(s)=\hat{F}_{B}(s)P_N\bigg(\int_0^{\infty}\hat{F}_{n-1}(ts)dF_T(t)\bigg),\ \ n\ge 1,
\end{eqnarray}
which is well-defined due to Lemma 2. Differentiating twice the
above equation with respect to $s$ and letting $s\to 0^+,$ we have,
for $n\ge 1,$
 \begin{eqnarray}
\hat{F}_n^{\prime}(0^+)&=&\hat{F}_B^{\prime}(0^+)+\hat{F}_{n-1}^{\prime}(0^+)\E[N]\E[T],\\
 \hat{F}_n^{\prime\prime}(0^+)&=&\hat{F}_B^{\prime\prime}(0^+)+2\hat{F}_B^{\prime}(0^+)\hat{F}_{n-1}^{\prime}(0^+)\E[N]\E[T]\nonumber\\
& &+\,\E[N(N-1)](\hat{F}_{n-1}^{\prime}(0^+)\E[T])^2+\hat{F}_{n-1}^{\prime\prime}(0^+)\E[N]\E[T^2].
\end{eqnarray}
\indent With the help of Lemma 6 and  by induction on $n,$ we can
show through (53) and (54) that $\E[Y_n]=\E[Y_0]=\mu_1,$
$\E[Y_n^2]=\E[Y_0^2]=\mu_2$ (defined in (50) and (51)) for all $n\ge
1$ and hence
\begin{eqnarray}\hbox{Var}(Y_n)=\mu_2-\mu_1^2=\frac{\hbox{Var}(B)+\mu^2(\E[T])^2 \hbox{Var}(N)+\mu^2\E[N]\hbox{Var}(T)}{1-\E[N]\E[T^2]}
,\ \ \ n\ge 0.
\end{eqnarray}
Moreover, by Lemma 4, we first have $\hat{F}_1\le \hat{F}_0$, and
then by the iteration (52), $\hat{F}_n\le \hat{F}_{n-1}$ for all
$n\ge 2.$ Namely, $\{Y_n\}_{n=0}^{\infty}$ is a sequence of
nonnegative random variables having the same first two moments
$\mu_1, \mu_2,$ and their Laplace--Stieltjes transforms
$\{\hat{F}_n\}$ are decreasing. Therefore, Lemma 9 applies. Denote
the limit of $\{\hat{F}_n\}$ by $\hat{F}_{\infty},$ which is the
Laplace--Stieltjes transform of a nonnegative random variable
$Y_{\infty}\sim F_{\infty}$ with $\E[Y_{\infty}]=\mu_1$ and
$\E[Y_{\infty}^2]\in[\mu_1^2,\mu_2].$ Consequently, it follows from
(52) that the limit $F_{\infty}$ is a solution to Eq (15), which has
mean
 $\mu$ and a
finite variance. Applying Lemma 6 to Eq (15) again, we conclude that
$\E[Y_{\infty}^2]=\mu_2$ as given in (51), and hence the solution
${Y}_{\infty}\sim F_{\infty}$ has the required variance as shown in
(17) or (55). This proves the necessity part of (a).

For part (b),  the proof of the uniqueness of the  solution to Eq
(15) is similar to that of Theorem 1,  and is omitted. The proof of
the theorem is complete.

\noindent{\bf Proof of Theorem 5.} In view of Eqs (1) and (18), we
want to prove the equivalence of the two functional equations:
\begin{eqnarray}
\hat{F}_*(s)=P_N\bigg(\int_0^{\infty}\hat{F}_*(ts)dF_T(t)\bigg),\ \ s\ge 0,
\end{eqnarray}
with $\lim_{s\to 0^+}(1-\hat{F}_{*}(s))/s=\mu\in(0,\infty),$ and
\begin{eqnarray}
\hat{F}_{\alpha}(s)=P_N\bigg(\int_0^{\infty}\hat{F}_{\alpha}(t^{1/{\alpha}}s)dF_T(t)\bigg),\ \ s\ge 0,
\end{eqnarray}
with $\lim_{s\to
0^+}(1-\hat{F}_{\alpha}(s))/s^{\alpha}=\mu\in(0,\infty).$ Here
\begin{eqnarray}
\hat{F}_{\alpha}(s)=\E\big[\exp({-sX_{\alpha}})\big]=\E\big[\exp({-sT_{\alpha}X_*^{1/\alpha}})\big]=\int_0^{\infty}e^{-s^{\alpha}x}dF_*(x)=\hat{F}_*(s^{\alpha}),\ \ s\ge 0.
\end{eqnarray}
\indent Suppose $\hat{F}_*$ satisfies (56) with $\lim_{s\to
0^+}(1-\hat{F}_{*}(s))/s=\mu.$ Then it follows from (58) that
\[\hat{F}_{\alpha}(s)=\hat{F}_*(s^{\alpha})=P_N\bigg(\int_0^{\infty}\hat{F}_*(ts^{\alpha})dF_T(t)\bigg)=P_N\bigg(\int_0^{\infty}\hat{F}_{\alpha}(t^{1/{\alpha}}s)dF_T(t)\bigg),\ \ s\ge 0,
\]
and $\lim_{s\to 0^+}(1-\hat{F}_{\alpha}(s))/s^{\alpha}=\lim_{s\to
0^+}(1-\hat{F}_*(s^{\alpha}))/s^{\alpha}=\lim_{s\to
0^+}(1-\hat{F}_{*}(s))/s=\mu.$ This means that (56) implies (57).
The converse implication can be proved similarly.
\bigskip\\
\noindent{\bf 5. Discussions}

 We have to mention that under the setting of Eq (1), the distributional equation is different from the following one
 (in which all i.i.d. $T_i$ are replaced by the same $T$):
\begin{eqnarray}
X\stackrel{\rm d}{=}T\sum_{i=1}^NX_i.
\end{eqnarray}
It is seen that Eq (59) has instead the corresponding functional
form
\begin{eqnarray}
\hat{F}(s)=\int_0^{\infty}P_N\big(\hat{F}(ts)\big)dF_T(t)=\int_0^{\infty}P_N\big(\E[\exp(-stX)]\big)dF_T(t),\ \ s\ge 0,
\end{eqnarray}
which is not equal to Eq (2) in general. But when $T$ is degenerate
at $p\in(0,1),$ Eq (60) also reduces to the Poincar\'e functional
equation (3) as Eq (2) does. Therefore, Eq (60) is another
generalization of the Poincar\'e functional equation (3).

The solutions to Eqs (2) and (60) are  distinct in general; in fact,
the second moments of the distributional solutions are different
from each other.  More precisely, the second moment of the solution
$X\sim F$ (with mean $\mu$) to Eq (60) is of the form
\[\E[X^2]=\frac{\E[N(N-1)]\E[T^2]}{1-\E[N]\E[T^2]}\cdot \mu^2,
\]
which is greater than or equal to that in (32) because $\E[T^2]\ge
(\E[T])^2.$

 It is interesting, however, that the necessary and sufficient
conditions for Eq (60) to have exactly one solution with finite
variance are the same as those for Eq (2). Namely, we have the
following result. The proof is similar to that of Theorem 1 and is
omitted.

\noindent{\bf Theorem 6.} Under the setting of Theorem 1 with given
$\mu,$  the random variables $N$ and $T$ together satisfy the
conditions (6) iff  the functional equation (60) has exactly one
solution $F$ with mean $\mu$ and a finite variance. Moreover, the
variance is of the form
\begin{eqnarray*}\hbox{Var}(X)=\frac{\E[N^2]\E[T^2]-1}{1-\E[N]\E[T^2]} \cdot \mu^2
\end{eqnarray*}
with $\E[N]=1/{\E[T]}.$
\medskip\\
\indent Theorem 5 about Eq (1) has a parallel result for Eq (59),
which extends both Theorems 1 and 2 of Hu and Cheng (2012). In this
regard, see also Hu and Lin (2001), Section 4, for characterizations
of the so-called semi-Mittag-Leffler distributions.

 Analogously, Eq (14) is
different from the following:
\begin{eqnarray*}
X\stackrel{\rm d}{=}B+T\sum_{i=1}^NX_i,
\end{eqnarray*}
which has the corresponding functional form
\begin{eqnarray}
& &\hat{F}(s)=\hat{F_B}(s) \cdot \int_0^{\infty}P_N\big(\hat{F}(ts)\big)dF_T(t)\\
&=&\hat{F_B}(s) \cdot \int_0^{\infty}P_N\big(\E[\exp(-stX)]\big)dF_T(t),\ \ s\ge 0.\nonumber
\end{eqnarray}
For this case, we have the next result analogous to  Theorem 4. The
proof is also omitted.

 \noindent{\bf Theorem 7.} Under the setting of Theorem 4,   the following statements are true.
\\
(a) The random variables $B, N$ and $T$ together satisfy the
conditions (16) iff the functional equation (61) has one solution
$F$ with mean $\mu$ and a finite variance. Moreover, the variance is
of the form
\begin{eqnarray*}
\hbox{Var}(X)=\frac{\hbox{Var}(B)+\E[B](2\mu-\E[B])+(\E[N^2]\E[T^2]-1)\mu^2}{1-\E[N]\E[T^2]}
\end{eqnarray*}
with  $\mu={\E[B]}/(1-\E[N]\E[T]).$
\\
(b) If, in addition to (16), $\E[T^2]<\E[T],$ then the solution $F$
to Eq (61) with given mean is unique.
\medskip\\
\indent
Finally, we remark that two different distributional equations may
be transferred to the same functional form which of course leads to
the same distributional solution. This means that a probability
distribution may have several kinds of characteristic properties.
\bigskip\\
 \centerline{\bf References}
\begin{description}

\item Azlarov, T.\,A. and Volodin, N.\,A. (1986). {\it Characterization Problems Associated with the Exponential Distribution.} Springer, New York.

\item  Eckberg, A.\,E.\,Jr. (1977). Sharp bounds on Laplace--Stieltjes transforms, with applications to
various queueing problems. {\it Math. Oper. Res.}, {\bf 2},
135--142.

\item Feller, W. (1971). {\it An Introduction to Probability Theory and its Applications}, Vol. 2, 2nd edn. Wiley,
New York.

\item Gulja$\check{\hbox{s}}$, B., Pearce,  C.\,E.\,M. and Pe$\check{\hbox{c}}$ari\'c, J. (1998). Jensen's inequality for distributions possessing
higher moments, with application to sharp bounds for
Laplace--Stieltjes transforms. {\it J. Austral. Math. Soc., Ser.
B}, {\bf 40}, 80--85.

\item Harkness, W.\,L. and Shantaram, R. (1969). Convergence of a
sequence of transformations of distribution functions. {\it
Pacific J. Math.}, {\bf 31}, 403--415.

\item Hu, C.-Y. and Cheng, T.-L. (2012). A characterization of distributions by random summation. {\it Taiwanese J. Math.}, {\bf 16}, 1245--1264.

\item Hu, C.-Y. and Lin, G.\,D. (2001). On the geometric compounding
model with applications. {\it Probab. Math. Statist.}, {\bf 21},
135--147.

\item Hu, C.-Y. and Lin, G.\,D. (2008). Some inequalities for
Laplace transforms. {\it J. Math. Anal. Appl.}, {\bf 340},
675--686.

\item Hu, C.-Y. and Lin, G.\,D. (2018). Characterizations of the logistic
and related distributions. {\it J. Math. Anal. Appl.}, {\bf
463}, 79--92.

\item Lin, G.\,D. (1993). Characterizations of the exponential distribution via the blocking time in a queueing system. {\it Statist. Sinica}, {\bf 3}, 577--581.

\item Lin, G.\,D. (1994). Characterizations of the Laplace and related distributions via geometric compound. {\it Sankhy$\bar{a}$, Ser. A}, {\bf 56}, 1--9.

\item Lin, G.\,D. (1998). Characterizations of the ${\cal L}$-class of life
distributions. {\it Statist. Probab. Lett.}, {\bf 40}, 259--266.

\item Lin, G.\,D. (2003). Characterizations of the exponential
distribution via the residual lifetime. {\it Sankhy$\bar{a}$,
Ser. A}, {\bf 65}, 249--258.

\item Liu, Q. (1997). Sur une \'equation fonctionnelle et ses applications: une extension du th\'eor\`eme de Kesten--Stigum concernant des processus de branchement.
[On a functional equation and its applications: an extension of
the Kesten--Stigum theorem on branching processes.] {\it Adv.
Appl. Probab.},  {\bf 29}, 353--373.

\item Liu, Q. (2002). An extension of a functional equation of Poincar\'e and Mandelbrot. {\it Asian J. Math.}, {\bf 6}, 145--168.

\item Poincar\'e, H. (1886). Sur une classe \'etendue de transcendantes uniformes. {\it C. R. Acad. Sci. Paris}, {\bf 103},
862--864.

\item Poincar\'e, H. (1890). Sur une classe nouvelle de transcendantes uniformes. {\it J. Math.
Pures Appl. IV. Ser.}, {\bf 6}, 313--365.

\item Ramachandran, B. and Lau, K.-S. (1991). {\it Functional Equations in Probability Theory.} Academic Press, New York.

\item Rao, C.\,R.,  Sapatinas, T. and Shanbhag, D.\,N. (1994). The integrated Cauchy functional equation: some comments on recent papers.
 {\it Adv. Appl. Probab.}, {\bf 26},  825--829.

 \item Steutel, F.\,W. and van Harn, K. (2004). {\it Infinite Divisibility of Probability Distributions on the Real Line.} Marcel Dekker, New York.

\end{description}
\end{document}